\documentclass[11pt,reqno]{amsart}
\usepackage{amsmath}
\usepackage{amsthm}
\usepackage{amssymb}
\usepackage{amsfonts}
\usepackage[foot]{amsaddr}
\usepackage{mathrsfs}
\usepackage{bm}
\usepackage{bbm}
\usepackage{color}
\usepackage{geometry}
\usepackage{graphicx}
\usepackage{hyperref}
\hypersetup{hypertex=true,
	colorlinks=true,
	linkcolor=blue,
	anchorcolor=blue,
	citecolor=blue}
\usepackage{aligned-overset}
\usepackage{yhmath}

\usepackage{autobreak}
\allowdisplaybreaks

\usepackage{cite}

\numberwithin{equation}{section}

\def\T{{\mathbb{T}}}

\def\R{{\mathbb{R}}}

\theoremstyle{definition}
\newtheorem{definition}{Definition}[section]

\newtheorem{theorem}[definition]{Theorem}
\theoremstyle{remark}
\newtheorem{remark}[definition]{Remark}

\begin{document}
\title[Self-dual solution of NS equations]{Self-dual solution of 3D incompressible Navier-Stokes equations}


\author{Ning-An Lai$^1$}
\address{$^1$School of Mathematical Sciences, Zhejiang Normal University, Jinhua 321004, China}

\author{Yi Zhou$^2$}
\address{$^2$School of Mathematical Sciences, Fudan University, Shanghai 200433, China. }

\email{ninganlai@zjnu.edu.cn, yizhou@fudan.edu.cn}

\date{\today}


\begin{abstract}

Whether the 3D incompressible Navier-Stokes equations will have a global smooth solution for all smooth, finite energy initial data is a Millennium Prize problem. One of the main difficulties of this problem is that the Navier-Stokes equations are actually a system of semilinear heat equations rather than a single equation. In this paper, we discover a remarkable hidden symmetry of the 3D incompressible Navier-Stokes equations. 
Under this symmetric reduction, the system reduces to a single scalar semilinear heat equation. The symmetry also holds for the 3D incompressible Euler equations.	

\end{abstract}

\maketitle

\section{Introduction}

The study of three-dimensional fluids remains an important issue nowadays. It has been considered as an important challenge for the 21st century. One of the major open problems is: Navier–Stokes Existence and Smoothness, which is one of the seven Millennium Prize Problems established by the Clay Mathematics Institute of Cambridge, Massachusetts (CMI) in 2000, in order to celebrate mathematics in the new millennium. The website describes it as:\\
\emph{Waves follow our boat as we meander across the lake, and turbulent air currents follow our flight in a modern jet. Mathematicians and physicists believe that an explanation for and the prediction of both the breeze and the turbulence can be found through an understanding of solutions to the Navier-Stokes equations. Although these equations were written down in the 19th Century, our understanding of them remains minimal. The challenge is to make substantial progress toward a mathematical theory which will unlock the secrets hidden in the Navier-Stokes equations.}

The official description of this Clay Millennium Problem was given by the 1978 Fields medalist Fefferman \cite{Fe06} as follows: we ask for a proof of one
of the following four statements\\
\begin{itemize}
		\item Existence and smoothness of Navier-Stokes solutions on $\R^3$;\\
        \item Existence and smoothness of Navier–Stokes solutions in $ \T^3$;\\
        \item Breakdown of Navier–Stokes solutions on $\R^3$;\\
        \item Breakdown of Navier–Stokes Solutions on $ \T^3$.
\end{itemize}
We also refer to the classic book \cite{Lem16} for the detailed description of this open problem.

In this paper we consider the Cauchy problem for the Navier-Stokes equations for an incompressible fluid:
\begin{equation}\label{ENS}
	\left\{
	\begin{aligned}
		& \partial_{t}u+u\cdot\nabla u+\nabla p =\nu \Delta u,~~~t\ge 0, x\in \mathbb{R}^3,\\
		&\nabla\cdot u=0,   ~~~~t\ge 0, x\in \mathbb{R}^3,\\
        & u(0, x)=u_0(x),~~~~~x\in \R^3,
	\end{aligned}
	\right.  
\end{equation}
where the unknown are the velocity vector $u(t, x)=(u_1(t, x), u_2(t, x), u_3(t, x))$ and the pressure $p(t, x)$, $\nu$ is a positive constant denoting the viscosity. Given the initial state $u_0(x)$, we want to determine the evolution of the system \eqref{ENS} for $t>0$. The local existence of classical solutions to the system \eqref{ENS} has been proved in \cite{Oseen11,Oseen27}, which was extended to the existence of global-in-time weak solutions by \cite{Le34}. Unfortunately, the global existence of classical solutions remains an important
challenge.

The known a priori Leray–Hopf energy estimate 
\begin{equation}\label{LHenergy}
\begin{aligned}
\frac12\int_{\mathbb{R}^3}|u(t, x)|^2dx+\nu\int_0^t\int_{\R^3}|\nabla u(s, x)|^2dxds=\frac12\int_{\R^3}|u(0, x)|^2dx
  \end{aligned}
\end{equation}
for the classical solutions to the system \eqref{ENS} can be obtained by a standard way. However, there is a "scaling gap" between this a priori bound and any known regularity criterion. Actually, if $(u, p)$ solves the system \eqref{ENS}, then also the scaled pair $(u^\lambda, p^\lambda)$
\[
u^\lambda(t, x)=\lambda u(\lambda t, \lambda x),~~ p^{\lambda}(t, x)=\lambda^2p
(\lambda t, \lambda x),~~~\lambda>0.
\]
As usual, we may assign the time $t$ a positive dimension 2, each space variable $x_i$ a positive dimension 1, $u$ a negative dimension $-1$ and $p$ a negative dimension -2. Then a direct dimensional analysis shows that all the terms in the a priori estimate \eqref{LHenergy} have positive dimension 1, which implies that the NS (short for Navier-Stokes) regularity problem is supercritical. 
As a result, any attempt to prove the global regularity for NS is fruitless at present. Thus, some of the research group has tried to search the blow up solution for the NS equations. The famous partial regularity result by Caffarelli, Kohn and Nirenberg \cite{CKN82} (which was simplified by Lin \cite{Lin98}) implies that any potential blow up of the axisymmetric Navier–
Stokes equations must occur on the symmetry axis. Another interesting result is due to Tao \cite{Tao16}, in which an averaged three-dimensional Navier–Stokes
equation preserving the energy identity but blowing up in finite time was proposed. We also refer to \cite{Hou23} for numerical analysis for potential singularity and \cite{Lem16} for a survey of the recent advances.

On the other hand, there exists a conserved quantity
\begin{equation}\label{Heli}
\begin{aligned}
H(u)=\frac12\int_{\R^3}u\cdot\omega dx+\nu\int_0^t\int_{\R^3}\nabla u\cdot\nabla\omega dxdt,\\
  \end{aligned}
\end{equation}
where $\omega=\nabla\times u$. $H(u)$ is the so called helicity satisfying 
\begin{equation}\label{Helicon}
\begin{aligned}
\frac{d H(u)}{dt}=0.
  \end{aligned}
\end{equation}
It is easy to see that $H(u)$ is critical with respect to the natural scaling of the Navier-Stokes equations. We believe that it is very important to study this quantity for the possible singularity or global existence to the  Navier-Stokes equations. In the recent paper by Lei, Lin and the second author \cite{LLZ15}, they discover a useful structure of the helicity, i.e. let 
\begin{equation}\label{u+u-}
\begin{aligned}
u_+=\frac12\left(u+\Lambda^{-1}\nabla\times u\right),\\
u_-=\frac12\left(u-\Lambda^{-1}\nabla\times u\right),
  \end{aligned}
\end{equation}
be the projection of the solutionn $u$ to the positive (respectively negative ) spectrum of the curl operator, then we have
\begin{equation}\label{prou+u-}
\begin{aligned}
&\nabla\times u_+=\Lambda u_+,\\
&\nabla\times u_-=-\Lambda u_-,\\
  \end{aligned}
\end{equation}
and
\begin{equation}\label{Ec1}
\begin{aligned}
&E_c(u(t))=E_c(u_+(t))+E_c(u_-(t)),\\
&H(u)=E_c(u_+(t))-E_c(u_-(t)),
  \end{aligned}
\end{equation}
where
\[
\Lambda=\sqrt{-\Delta}
\]
and
\begin{equation}\label{Ec}
\begin{aligned}
E_c(u(t)):=\frac12\left\|\Lambda^{\frac12}u(t, \cdot)\right\|_{L^2}^2+
\nu\int_0^t\left\|\Lambda^\frac{1}{2}\nabla u(s, \cdot)\right\|_{L^2}^2ds
  \end{aligned}
\end{equation}
denotes the critical energy, which is invariant under the natural scaling of the Navier-Stokes equations. In \cite{LLZ15}, they found that in the extremal case 
\[
E_c(u_+)\equiv 0~ or ~E_c(u_-)\equiv 0, 
\]
then $u_-$ or $u_+$ becomes the Baltrami flow, i.e.
\[
\nabla\times u_{\pm}=\pm k u_{\pm},
\] 
where $k$ is a positive real constant and the Navier-Stokes equations become a linear system, hence the global existence follows, see \cite{CM88,MB02}. In \cite{LLZ15, LZZ16arxiv, WZZ24}, a various sorts of initial data were found such that
\begin{equation}\label{Data}
\begin{aligned}
E_c(u_+)(t)\ll E_c(u_-)(t),~~~\forall t\ge 0,\\
  \end{aligned}
\end{equation}
or 
\begin{equation}\label{Data2}
\begin{aligned}
E_c(u_-)(t)\ll E_c(u_+)(t),~~~\forall t\ge 0,\\
  \end{aligned}
\end{equation}
and it has been proved there is always global existence of smooth solution in this case.
By common sense, it is better to search for the possible blow up solution in the other extremal case $H(u)\equiv 0$. In this paper, we discover a remarkable hidden symmetry of the 3D Navier-Stokes equation,  i.e. under the transformation
\[
\begin{aligned}
&u_+(t, x)\rightarrow -u_-(t, -x),\\
&u_-(t, x)\rightarrow -u_+(t, -x),\\
\end{aligned}
\] 
the 3D Navier-Stokes equation is invariant, and this allows us to define what we call the self-dual solution
\begin{equation}
u_-(t,x)=-u_+(t,-x),
\end{equation}
which implies
\[
E_c(u_+)(t)=E_c(u_-)(t),
\]
and hence $H(u)=0$. 

By \eqref{prou+u-}, we see that the space Fourier transform of $u_+$ has only one degree of freedom, so we can write 
\begin{equation}
u_+(t,x)=v(t,x)\ast h(x),
\end{equation}
where $v$ is a real scalar function of $(t,x)$ and $h$ is any real distribution from $\R^3$ to $\R^3$ satisfying
\begin{equation}\label{h}
\begin{aligned}
&\nabla\times h(x)=\Lambda h(x),\\
&\nabla\cdot h(x)=0.
\end{aligned}
\end{equation}
If $h$ satisfies further
\begin{equation}
|\widehat{h}|=1,
\end{equation}
where $\widehat{h}$ denotes the Fourier transform of $h$, then we will show in the next section that $v$ satisfies
\begin{equation}\label{veq}
	\left\{
	\begin{aligned}
		& v_t-\nu \Delta v=\sum_{i=1}^3h_i^r\ast\left[\left(v\ast h-v^r\ast h^r\right)\times
\left(v\ast \Lambda h+v^r\ast\Lambda h^r\right)\right]_i,~~~t\ge 0, x\in \mathbb{R}^3,\\
		&v^r(t, x)=v(t, -x), h^r(t, x)=h(t, -x),   ~~~~t\ge 0, x\in \mathbb{R}^3.\\
	\end{aligned}
	\right.  
\end{equation}



Thus, we reduce the original nonlinear system to a single scalar semilinear heat equation. Compared to the nonlinear system of 3D Naver-Stokes equations, it will be much easier to study a single semilinear heat equation, and hence this reduction may shed some new light to the 
Millennium Prize problem. Moreover, by rotational invariance of the equation, we may consider an axisymmetric solution which is even simpler to investigate. By the reduction in the next section, we arrive at the following theorem
\begin{theorem}\label{thm1}
Let
\[
\Lambda'=\sqrt{-\partial_r^2-\frac{\partial_r}{r}}.
\]
Suppose that there is smooth scalar function $v_0(r, x_3)\in L^2(\R^3)$ which is axial symmetric such that the scalar equation 
\begin{equation}\label{nsmain}
\begin{aligned}
&2\sqrt 2\left[v_t-\nu\left(\partial_r^2+\frac{\partial_r}{r}+\partial_{x_3}^2\right)
v\right]\\
=&\Lambda'\Lambda^{-1}\left\{\left[(\Lambda')^{-1}(v+v^r)
\right]_r\left[(\Lambda')
^{-1}(v+v^r)\right]_{rx_3}\right\}\\
&-\Lambda'\Lambda^{-1}\left\{\left[(\Lambda')
^{-1}\Lambda^{-1}(v-v^r)
\right]_{rx_3}\left[\Lambda(\Lambda')
^{-1}(v-v^r)\right]_r\right\}\\
&+(\Lambda')^{-1}\Lambda^{-1}\partial_{x_3}\left\{\left[
\Lambda^{-1}\Lambda'(v-v^r)\right]_r\left[\Lambda(\Lambda')^{-1}(v-v^r)\right]_r
\right\}\\
&-(\Lambda')^{-1}\Lambda^{-1}\partial_{x_3}\left[\Lambda^{-1}\Lambda'(v-v^r)
\Lambda\Lambda'(v-v^r)\right]-(\Lambda')^{-1}\left\{\left[\Lambda'
\Lambda^{-1}(v-v^r)\right]_r
\left[(\Lambda')^{-1}(v+v^r)\right]_{rx_3}\right\}\\
&+(\Lambda')^{-1}\left\{\Lambda'\Lambda^{-1}(v-v^r)\left[\Lambda'(v+v^r)
\right]_{x_3}\right\}\\
&-\partial_{x_3}(\Lambda')^{-1}\Lambda^{-1}\left\{\left[(\Lambda')^{-1}(v+v^r)\right]_r
\left[\Lambda'(v-v^r)\right]_r\right\}\\
&+\partial_{x_3}(\Lambda')^{-1}\Lambda^{-1}\left\{\Lambda'(v+v^r)
\Lambda'(v+v^r)\right\}\\
&+(\Lambda')^{-1}\left\{\left[(\Lambda')^{-1}\Lambda^{-1}(v-v^r)\right]_{rx_3}
\left[\Lambda'(v+v^r)\right]_r\right\}\\
&-(\Lambda')^{-1}\left\{\left[\Lambda'\Lambda^{-1}(v-v^r)\right]_{x_3}
\Lambda'(v+v^r)\right\}
  \end{aligned}
\end{equation}
with initial data $v_0$ blows up in a finite time, then we can find a smooth divergence free vector field $u_0(x)\in L^2(\R^3)$ such that the Cauchy problem \eqref{ENS} of the Navier-Stokes will also blow up in a finite time.

\end{theorem}


\begin{remark} 

If $\nu=0$, then the system \eqref{ENS} becomes the incompressible Euler system, the hidden symmetry also holds for this case, and hence Theorem \ref{thm1} also holds for the 3D incompressible Euler system. The global regularity for the 3D incompressible Euler equations is also one of the most outstanding open problem in the past century, see the survey paper \cite{BT07, Cha08, Con07} and lecture notes \cite{HY}.
\end{remark}

For the singularity formation analysis, it is very important to study the corresponding stationary solution. If a stationary solution with appropriate regularity to \eqref{nsmain} can be found, then we may possibly use the modulation method to construct the finite time blow up solution to the Navier-Stokes equations, see the details in Section 3.

\section{The reduced single equation}

Let 
\[
\psi(x): \R^3\rightarrow \R^3
\]
be a flow map and
\[
Q\psi=(-\Delta)^{-1}\nabla\times\nabla\times \psi
\]
be the projection to the divergence free field. For $\phi$ satisfying 
\[
\nabla\cdot\phi=0,
\]
define
\[
Q_{\pm}\phi=\frac12\left(\phi\pm\Lambda^{-1}\nabla\times\phi\right)
\]
and let
\[
P_{\pm}=Q_{\pm}Q.
\]
Since the Navier-Stokes equations can be written as
\begin{equation}\label{ns}
\begin{aligned}
u_t-\nu\Delta u-u\times\omega+\nabla\left(\frac{|u|^2}{2}+p\right)=0,
  \end{aligned}
\end{equation}
then by \eqref{prou+u-} we have
\begin{equation}\label{ns1}
\begin{aligned}
u_t-\nu\Delta u+\nabla \left(\frac{|u|^2}{2}+p\right)=(u_++u_-)\times(\Lambda u_+-\Lambda u_-),
  \end{aligned}
\end{equation}
which yields
\begin{equation}\label{ns2}
\begin{aligned}
&\partial_tu_{+}-\nu\Delta u_+=P_+\left[(u_++u_-)\times(\Lambda u_+-\Lambda u_-)\right],\\
&\partial_tu_{-}-\nu\Delta u_-=P_-\left[(u_++u_-)\times(\Lambda u_+-\Lambda u_-)\right].\\
  \end{aligned}
\end{equation}
Under the transformation
\[
u_+\rightarrow -u_-^r, u_-\rightarrow-u_+^r,
\]
the equations \eqref{ns2} are invariant. So we can consider a self-dual solution
\[
u_-=-u_+^r,
\]
and then the equation $\eqref{ns2}_1$ becomes 
\begin{equation}\label{ns5}
\begin{aligned}
&\partial_tu_{+}-\nu\Delta u_+=P_+\left[(u_+-u_+^r)\times(\Lambda u_++\Lambda u_+^r)\right].\\
  \end{aligned}
\end{equation}

Since
\[
\begin{aligned}
&\nabla\times u_+=\Lambda u_+,\\
&\nabla\cdot u_+=0,
  \end{aligned}
\]
which yields
\[
\begin{aligned}
&\sqrt{-1}\xi\times \widehat{u}_+=|\xi|\widehat{u}_+,\\
&\xi \cdot \widehat{u}_+=0.
  \end{aligned}
\]
So $\widehat{u}_+(\xi)$ is the eigenvector of
\begin{equation*}
A(\xi)=\sqrt{-1}\begin{pmatrix}
0 & -\xi_3 & +\xi_2\\
\xi_3 & 0 & -\xi_1\\
-\xi_2 &+\xi_1 &0\\
\end{pmatrix}
\end{equation*}
with respect to eigenvalue $|\xi|$ and then has only one degree of freedom. We then may write 
\[
\widehat{u}_+(t, \xi)=\widehat{v}(t, \xi)\widehat{h}(\xi),
\]
where $h(x)$ satisfies \eqref{h}. And hence
\begin{equation}\label{uh}
u_+=v\ast h.
\end{equation}
\begin{equation}\label{ns5}
\begin{aligned}
&\sum_{i}\left(\partial_tu_+-\nu\Delta u_+\right)_i\ast h_i^r=
\sum_ih_i^r\ast\left(P_+\left[(u_+-u_+^r)\times(\Lambda u_++\Lambda u_+^r)\right]\right)_i\\
&=(P_+h)^r_i\ast\left[(u_+-u_+^r)\times(\Lambda u_++\Lambda u_+^r)\right]_i,\\
  \end{aligned}
\end{equation}
since
\[
P_+h=h,
\]
so the right hand side of the above equality is
\begin{equation}\label{ns5a}
\begin{aligned}
\sum_{i}h^r_i\ast\left[(u_+-u_+^r)\times(\Lambda u_++\Lambda u_+^r)\right]_i.\\
  \end{aligned}
\end{equation}
The Fourier transform of the left hand side of \eqref{ns5} equals to
\begin{equation}\label{ns6}
\begin{aligned}
&\left(\partial_t\widehat{u}_++\nu|\xi|^2\widehat{u}_+\right)\cdot \overline{\widehat{h}}\\
&=\left(\widehat{v}_t+\nu|\xi|^2\widehat{v}\right)\widehat{h}
\cdot\overline{\widehat{h}}\\
&=\widehat{v}_t+\nu|\xi|^2\widehat{v},
  \end{aligned}
\end{equation}
which leads to 
\begin{equation}\label{lev}
v_t-\nu\Delta v
\end{equation}
after making inverse Fourier transform. Inserting \eqref{uh} into \eqref{ns5}, and combining \eqref{lev}, we get the equation \eqref{veq}.

Making Fourier transform with respect to space variable to the equation \eqref{veq}, we get
\begin{equation}\label{ns6}
\begin{aligned}
\partial_t\widehat{v}(t, \xi)+\nu|\xi|^2\widehat{v}(t, \xi)=&\int_{\R^3}|\eta|\overline{\widehat{h}}(\xi)\cdot\left[\widehat{h}(\xi-\eta)\times \widehat{h}(\eta)\right]
\widehat{v}(t, \xi-\eta)\widehat{v}(\eta)d\eta\\
&-\int_{\R^3}|\eta|\overline{\widehat{h}}(\xi)\cdot\left[\overline{\widehat{h}}
(\xi-\eta)\times \widehat{h}(\eta)\right]
\overline{\widehat{v}}(t, \xi-\eta)\widehat{v}(\eta)d\eta\\
&+\int_{\R^3}|\eta|\overline{\widehat{h}}(\xi)\cdot\left[\widehat{h}(\xi-\eta)\times \overline{\widehat{h}}(\eta)\right]
\widehat{v}(t, \xi-\eta)\overline{\widehat{v}}(\eta)d\eta\\
&-\int_{\R^3}|\eta|\overline{\widehat{h}}(\xi)\cdot\left[\overline{\widehat{h}}
(\xi-\eta)\times \overline{\widehat{h}}(\eta)\right]
\overline{\widehat{v}}(t, \xi-\eta)\overline{\widehat{v}}(\eta)d\eta.\\
  \end{aligned}
\end{equation}

Set
\begin{equation}\label{xi}
\begin{aligned}
&d(\xi)=(-\xi_2, \xi_1, 0),\\
&g(\xi)=\sqrt{-1}d(\xi)-|\xi|^{-1}\xi\times d(\xi),\\
&\widehat{h}(\xi)=\frac{g(\xi)}{|g(\xi)|}.\\
 \end{aligned}
\end{equation}
Direct computations yields 
\begin{equation}\label{coe}
\begin{aligned}
&\overline{\widehat{h}}(\xi)\left[\widehat{h}(\zeta)\times \widehat{h}(\eta)\right]\\
=&\frac{1}{2\sqrt{2}|\xi'||\zeta'||\eta'|}
\Bigg\{\frac{|\xi'|^2(\zeta_1\eta_2-\zeta_2\eta_1)}{|\xi|}\left(1-
\frac{\zeta_3\eta_3}{|\zeta||\eta|}\right)+\frac{|\zeta'|^2(\xi_2\eta_1
-\xi_1\eta_2)}{|\zeta|}\left(1-
\frac{\xi_3\eta_3}{|\xi||\eta|}\right)\\
&+\frac{|\eta'|^2(\xi_1\zeta_2-\xi_2\zeta_1)}{|\eta|}\left(1-
\frac{\xi_3\zeta_3}{|\xi||\zeta|}\right)\\
&-i\Bigg[(|\xi'|^2)(\zeta'\cdot\eta')\left(\frac{\eta_3}{|\xi||\eta|}
-\frac{\zeta_3}{|\xi||\zeta|}\right)+(|\zeta'|^2)(\xi'\cdot\eta')
\left(\frac{\xi_3}{|\xi||\zeta|}
-\frac{\eta_3}{|\zeta||\eta|}\right)\\
&-(|\eta'|^2)(\xi'\cdot\zeta')\left(\frac{\xi_3}{|\xi||\eta|}
-\frac{\zeta_3}{|\zeta||\eta|}\right)\Bigg]\Bigg\},
 \end{aligned}
\end{equation}
\begin{equation}\label{coe2}
\begin{aligned}
&\overline{\widehat{h}}(\xi)\left[\overline{\widehat{h}}(\zeta)\times \widehat{h}(\eta)\right]\\
=&\frac{1}{2\sqrt{2}|\xi'||\zeta'||\eta'|}
\Bigg\{\frac{|\xi'|^2(-\zeta_1\eta_2+\zeta_2\eta_1)}{|\xi|}\left(1+
\frac{\zeta_3\eta_3}{|\zeta||\eta|}\right)+\frac{|\zeta'|^2(\xi_2\eta_1
-\xi_1\eta_2)}{|\zeta|}\left(1-
\frac{\xi_3\eta_3}{|\xi||\eta|}\right)\\
&+\frac{|\eta'|^2(-\xi_1\zeta_2+\xi_2\zeta_1)}{|\eta|}\left(1+
\frac{\xi_3\zeta_3}{|\xi||\zeta|}\right)\\
&-i\Bigg[-(|\xi'|^2)(\zeta'\cdot\eta')\left(\frac{\eta_3}{|\xi||\eta|}
+\frac{\zeta_3}{|\xi||\zeta|}\right)+(|\zeta'|^2)(\xi'\cdot\eta')
\left(\frac{\xi_3}{|\xi||\zeta|}
-\frac{\eta_3}{|\zeta||\eta|}\right)\\
&+(|\eta'|^2)(\xi'\cdot\zeta')\left(\frac{\xi_3}{|\xi||\eta|}
+\frac{\zeta_3}{|\zeta||\eta|}\right)\Bigg]\Bigg\},
 \end{aligned}
\end{equation}
\begin{equation}\label{coe3}
\begin{aligned}
&\overline{\widehat{h}}(\xi)\left[\widehat{h}(\zeta)\times \overline{\widehat{h}}(\eta)\right]\\
=&\frac{1}{2\sqrt{2}|\xi'||\zeta'||\eta'|}
\Bigg\{\frac{|\xi'|^2(-\zeta_1\eta_2+\zeta_2\eta_1)}{|\xi|}\left(1+
\frac{\zeta_3\eta_3}{|\zeta||\eta|}\right)+\frac{|\zeta'|^2(-\xi_2\eta_1
+\xi_1\eta_2)}{|\zeta|}\left(1+
\frac{\xi_3\eta_3}{|\xi||\eta|}\right)\\
&+\frac{|\eta'|^2(-\xi_1\zeta_2+\xi_2\zeta_1)}{|\eta|}\left(-1+
\frac{\xi_3\zeta_3}{|\xi||\zeta|}\right)\\
&-i\Bigg[(|\xi'|^2)(\zeta'\cdot\eta')\left(\frac{\eta_3}{|\xi||\eta|}
+\frac{\zeta_3}{|\xi||\zeta|}\right)-(|\zeta'|^2)(\xi'\cdot\eta')
\left(\frac{\xi_3}{|\xi||\zeta|}
+\frac{\eta_3}{|\zeta||\eta|}\right)\\
&+(|\eta'|^2)(\xi'\cdot\zeta')\left(-\frac{\xi_3}{|\xi||\eta|}
+\frac{\zeta_3}{|\zeta||\eta|}\right)\Bigg]\Bigg\},
 \end{aligned}
\end{equation}
\begin{equation}\label{coe4}
\begin{aligned}
&\overline{\widehat{h}}(\xi)\left[\overline{\widehat{h}}(\zeta)\times \overline{\widehat{h}}(\eta)\right]\\
=&\frac{1}{2\sqrt{2}|\xi'||\zeta'||\eta'|}
\Bigg\{\frac{|\xi'|^2(\zeta_1\eta_2-\zeta_2\eta_1)}{|\xi|}\left(1-
\frac{\zeta_3\eta_3}{|\zeta||\eta|}\right)+\frac{|\zeta'|^2(-\xi_2\eta_1
+\xi_1\eta_2)}{|\zeta|}\left(1+
\frac{\xi_3\eta_3}{|\xi||\eta|}\right)\\
&+\frac{|\eta'|^2(-\xi_1\zeta_2+\xi_2\zeta_1)}{|\eta|}\left(1+
\frac{\xi_3\zeta_3}{|\xi||\zeta|}\right)\\
&-i\Bigg[-(|\xi'|^2)(\zeta'\cdot\eta')\left(\frac{\eta_3}{|\xi||\eta|}
-\frac{\zeta_3}{|\xi||\zeta|}\right)-(|\zeta'|^2)(\xi'\cdot\eta')
\left(\frac{\xi_3}{|\xi||\zeta|}
+\frac{\eta_3}{|\zeta||\eta|}\right)\\
&+(|\eta'|^2)(\xi'\cdot\zeta')\left(\frac{\xi_3}{|\xi||\eta|}
+\frac{\zeta_3}{|\zeta||\eta|}\right)\Bigg]\Bigg\}.
 \end{aligned}
\end{equation}
If we consider the axial symmetric solution 
\begin{equation}\label{axisol}
v(t, x)=v(t, \sqrt{x_1^2+x_2^2}, x_3),
\end{equation}
the Fourier transform of which can be written as
\[
\widehat{v}(t, \xi)=c(t, |\xi'|, \xi_3),
\]
where 
\[
\xi'=(\xi_1, \xi_2),
\]
then by making the inverse Fourier transform to \eqref{ns6}, we get
\begin{equation}\label{ns7}
\begin{aligned}
&v_t-\nu\left(\partial_r^2+\frac{\partial_r}{r}+\partial_{x_3}^2\right)v\\
=&\frac{1}{(2\pi)^3}\int_{\R^3}\int_{\R^3}e^{i\left[r(\zeta_1+\eta_1)+x_3(\zeta_3+
\eta_3)
\right]}
|\eta|\overline{\widehat{h}}(\zeta+\eta)\cdot
\left[\widehat{h}(\zeta)\times \widehat{h}(\eta)\right]
c(t, |\zeta'|, \zeta_3)c(t, |\eta'|, \eta_3)d\eta d\zeta\\
&-\frac{1}{(2\pi)^3}\int_{\R^3}\int_{\R^3}e^{i\left[r(\zeta_1+\eta_1)+x_3(\zeta_3+
\eta_3)
\right]}
|\eta|\overline{\widehat{h}}(\zeta+\eta)\cdot
\left[\overline{\widehat{h}}(\zeta)\times \widehat{h}(\eta)\right]
c(t, |\zeta'|, -\zeta_3)c(t, |\eta'|, \eta_3)d\eta d\zeta\\
&+\frac{1}{(2\pi)^3}\int_{\R^3}\int_{\R^3}e^{i\left[r(\zeta_1+\eta_1)+x_3(\zeta_3+
\eta_3)
\right]}
|\eta|\overline{\widehat{h}}(\zeta+\eta)\cdot
\left[\widehat{h}(\zeta)\times \overline{\widehat{h}}(\eta)\right]
c(t, |\zeta'|, \zeta_3)c(t, |\eta'|, -\eta_3)d\eta d\zeta\\
&-\frac{1}{(2\pi)^3}\int_{\R^3}\int_{\R^3}e^{i\left[r(\zeta_1+\eta_1)+x_3(\zeta_3+
\eta_3)
\right]}
|\eta|\overline{\widehat{h}}(\zeta+\eta)\cdot
\left[\overline{\widehat{h}}(\zeta)\times \overline{\widehat{h}}(\eta)\right]
c(t, |\zeta'|, -\zeta_3)c(t, |\eta'|, -\eta_3)d\eta d\zeta.\\
  \end{aligned}
\end{equation}
 Under the transformation 
\[
\begin{aligned}
&\zeta\rightarrow (\zeta_1, -\zeta_2, \zeta_3),\\
&\eta\rightarrow (\eta_1, -\eta_2, \eta_3),\\
\end{aligned}
\]
if a coefficient in the right hand side of \eqref{ns7} changes sign, then the corresponding term must vanish. In particular, equation \eqref{ns7} becomes
\begin{equation}\label{ns8}
\begin{aligned}
&v_t-\nu\left(\partial_r^2+\frac{\partial_r}{r}+\partial_{x_3}^2\right)v\\
=&\frac{\sqrt{-1}}{(2\pi)^3}\int_{\R^3}\int_{\R^3}e^{\sqrt{-1}\left[r(\zeta_1+\eta_1)
+x_3(\zeta_3+
\eta_3)
\right]}
Im\left\{|\eta|\overline{\widehat{h}}(\zeta+\eta)\cdot
\left[\widehat{h}(\zeta)\times \widehat{h}(\eta)\right]\right\}
c(t, |\zeta'|, \zeta_3)c(t, |\eta'|, \eta_3)d\eta d\zeta\\
&-\frac{\sqrt{-1}}{(2\pi)^3}\int_{\R^3}\int_{\R^3}e^{\sqrt{-1}\left[r(\zeta_1+\eta_1)+x_3(\zeta_3+
\eta_3)
\right]}
Im\left\{|\eta|\overline{\widehat{h}}(\zeta+\eta)\cdot
\left[\overline{\widehat{h}}(\zeta)\times \widehat{h}(\eta)\right]\right\}
c(t, |\zeta'|, -\zeta_3)c(t, |\eta'|, \eta_3)d\eta d\zeta\\
&+\frac{\sqrt{-1}}{(2\pi)^3}\int_{\R^3}\int_{\R^3}e^{\sqrt{-1}\left[r(\zeta_1+\eta_1)+x_3(\zeta_3+
\eta_3)
\right]}
Im\left\{|\eta|\overline{\widehat{h}}(\zeta+\eta)\cdot
\left[\widehat{h}(\zeta)\times \overline{\widehat{h}}(\eta)\right]\right\}
c(t, |\zeta'|, \zeta_3)c(t, |\eta'|, -\eta_3)d\eta d\zeta\\
&-\frac{\sqrt{-1}}{(2\pi)^3}\int_{\R^3}\int_{\R^3}e^{\sqrt{-1}\left[r(\zeta_1+\eta_1)+x_3(\zeta_3+
\eta_3)
\right]}
Im\left\{|\eta|\overline{\widehat{h}}(\zeta+\eta)\cdot
\left[\overline{\widehat{h}}(\zeta)\times \overline{\widehat{h}}(\eta)\right]\right\}
c(t, |\zeta'|, -\zeta_3)c(t, |\eta'|, -\eta_3)d\eta d\zeta,\\
  \end{aligned}
\end{equation}
which is still very complicated. However, we can write it more explicitly
\begin{equation}\label{ns9}
\begin{aligned}
&2\sqrt 2\left[v_t-\nu\left(\partial_r^2+\frac{\partial_r}{r}+\partial_{x_3}^2\right)
v\right]\\
=&\Lambda'\Lambda^{-1}\left\{\left[(\Lambda')^{-1}(v+v^r)
\right]_r\left[(\Lambda')
^{-1}(v+v^r)\right]_{rx_3}\right\}\\
&-\Lambda'\Lambda^{-1}\left\{\left[(\Lambda')
^{-1}\Lambda^{-1}(v-v^r)
\right]_{rx_3}\left[\Lambda(\Lambda')
^{-1}(v-v^r)\right]_r\right\}\\
&+(\Lambda')^{-1}\Lambda^{-1}\partial_{x_3}\left\{\left[
\Lambda^{-1}\Lambda'(v-v^r)\right]_r\left[\Lambda(\Lambda')^{-1}(v-v^r)\right]_r
\right\}\\
&-(\Lambda')^{-1}\Lambda^{-1}\partial_{x_3}\left[\Lambda^{-1}\Lambda'(v-v^r)
\Lambda\Lambda'(v-v^r)\right]-(\Lambda')^{-1}\left\{\left[\Lambda'
\Lambda^{-1}(v-v^r)\right]_r
\left[(\Lambda')^{-1}(v+v^r)\right]_{rx_3}\right\}\\
&+(\Lambda')^{-1}\left\{\Lambda'\Lambda^{-1}(v-v^r)\left[\Lambda'(v+v^r)
\right]_{x_3}\right\}\\
&-\partial_{x_3}(\Lambda')^{-1}\Lambda^{-1}\left\{\left[(\Lambda')^{-1}(v+v^r)\right]_r
\left[\Lambda'(v-v^r)\right]_r\right\}\\
&+\partial_{x_3}(\Lambda')^{-1}\Lambda^{-1}\left\{\Lambda'(v+v^r)
\Lambda'(v+v^r)\right\}\\
&+(\Lambda')^{-1}\left\{\left[(\Lambda')^{-1}\Lambda^{-1}(v-v^r)\right]_{rx_3}
\left[\Lambda'(v+v^r)\right]_r\right\}\\
&-(\Lambda')^{-1}\left\{\left[\Lambda'\Lambda^{-1}(v-v^r)\right]_{x_3}
\Lambda'(v+v^r)\right\},
  \end{aligned}
\end{equation}
which is exactly the one \eqref{nsmain}. We may
separate the odd and even parts to get
\begin{equation}\label{ns10}
\begin{aligned}
&\sqrt 2\left(\partial_t-\nu\Delta\right)(v-v^r)\\
=&\Lambda'\Lambda^{-1}\left\{\left[(\Lambda')^{-1}(v+v^r)
\right]_r\left[(\Lambda')
^{-1}(v+v^r)\right]_{rx_3}\right\}\\
&-\Lambda'\Lambda^{-1}\left\{\left[(\Lambda')
^{-1}\Lambda^{-1}(v-v^r)
\right]_{rx_3}\left[\Lambda(\Lambda')
^{-1}(v-v^r)\right]_r\right\}\\
&+(\Lambda')^{-1}\Lambda^{-1}\partial_{x_3}\left\{\left[
\Lambda^{-1}\Lambda'(v-v^r)\right]_r\left[\Lambda(\Lambda')^{-1}(v-v^r)\right]_r
\right\}\\
&-(\Lambda')^{-1}\Lambda^{-1}\partial_{x_3}\left[\Lambda^{-1}\Lambda'(v-v^r)
\Lambda\Lambda'(v-v^r)\right]\\
&+\partial_{x_3}(\Lambda')^{-1}\Lambda^{-1}\left\{\Lambda'(v+v^r)
\Lambda'(v+v^r)\right\}\\
  \end{aligned}
\end{equation}
and
\begin{equation}\label{ns11}
\begin{aligned}
&\sqrt 2\left(\partial_t-\nu\Delta\right)(v+v^r)\\
=&-(\Lambda')^{-1}\left\{\left[\Lambda'
\Lambda^{-1}(v-v^r)\right]_r
\left[(\Lambda')^{-1}(v+v^r)\right]_{rx_3}\right\}\\
&+(\Lambda')^{-1}\left\{\Lambda'\Lambda^{-1}(v-v^r)\left[\Lambda'(v+v^r)
\right]_{x_3}\right\}\\
&-\partial_{x_3}(\Lambda')^{-1}\Lambda^{-1}\left\{\left[(\Lambda')^{-1}(v+v^r)\right]_r
\left[\Lambda'(v-v^r)\right]_r\right\}\\
&+(\Lambda')^{-1}\left\{\left[(\Lambda')^{-1}\Lambda^{-1}(v-v^r)\right]_{rx_3}
\left[\Lambda'(v+v^r)\right]_r\right\}\\
&-(\Lambda')^{-1}\left\{\left[\Lambda'\Lambda^{-1}(v-v^r)\right]_{x_3}
\Lambda'(v+v^r)\right\}.
 \end{aligned}
\end{equation}
\begin{remark}
If we set 
\[
v+v^r=0,
\]
which is a solution to \eqref{ns11}, then the equation \eqref{ns10} for $V=v-v^r$ will change to
\begin{equation}\label{ns12}
\begin{aligned}
&\left(\partial_t-\nu\Delta\right)V\\
=&\frac{\sqrt2}{2}\Bigg(
-\Lambda'\Lambda^{-1}\left\{\left[(\Lambda')
^{-1}\Lambda^{-1}V
\right]_{rx_3}\left[\Lambda(\Lambda')
^{-1}V\right]_r\right\}\\
&+(\Lambda')^{-1}\Lambda^{-1}\partial_{x_3}\left\{\left[
\Lambda^{-1}\Lambda' V\right]_r\left[\Lambda(\Lambda')^{-1}V\right]_r
\right\}\\
&-(\Lambda')^{-1}\Lambda^{-1}\partial_{x_3}\left[\Lambda^{-1}\Lambda'V
\Lambda\Lambda'V\right]\Bigg).\\
  \end{aligned}
\end{equation}
The above semilinear heat equation looks much simpler, unfortunately, the solution to the corresponding Navier-Stokes system is swirl free and hence will not blow up.
\end{remark}

\section{Stationary solution}

The key ingredient to prove possible singularity is to find an appropriate stationary solution of the reduced single equation \eqref{ns9}. Set
\[
(\Lambda')^{-1}\Lambda^{-1}v=w,~~(\Lambda')^{-1}\Lambda^{-1}v^r=w^r,
\]
then by \eqref{ns9} the stationary equation for $w$ admits the form
\begin{equation}\label{ns13}
\begin{aligned}
&-2\sqrt2\nu\Delta\Delta\Delta' w\\
=&-\Delta'\left\{\left[\Lambda(w+w^r)\right]_r\left[\Lambda(w+w^r)\right]_{rx_3}
\right\}-\Delta'\left\{
\left(w-w^r\right)_{rx_3}\left[\Delta(w-w^r)\right]_r\right\}\\
&+\partial_{x_3}\left\{\left[\Delta'(w-w^r)\right]_r\left[\Delta(w-w^r)
\right]_r\right\}+\partial_{x_3}\left[\Delta'(w-w^r)\Delta\Delta'(w-w^r)\right]\\
&+\Lambda\left\{\left[\Delta'(w-w^r)\right]_r\left[\Lambda(w+w^r)
\right]_{rx_3}\right\}+\Lambda\left\{\Delta'(w-w^r)\left[\Lambda\Delta'(w+w^r)
\right]_{x_3}\right\}\\
&+\partial_{x_3}\left\{\left[\Lambda(w+w^r)\right]_r
\left[\Lambda\Delta'(w-w^r)
\right]_r\right\}+\partial_{x_3}\left[\Lambda\Delta'(w+w^r)\right]^2\\
&-\Lambda\left\{(w-w^r)_{rx_3}
\left[\Lambda\Delta'(w+w^r)
\right]_r\right\}-\Lambda\left\{\left[\Delta'(w-w^r)\right]_{x_3}
\Lambda\Delta'(w+w^r)\right\}.\\
  \end{aligned}
\end{equation}
We hope to find some nontrivial solution of \eqref{ns13} with bounded derivative $w_r$, and then the $BMO^{-1}$ norm of the solution $v$ can be controlled. Indeed,
\begin{equation}\label{ns17}
\begin{aligned}
\|v\|_{BMO^{-1}}&=\left\|\Lambda'w\right\|_{BMO}\\
&\lesssim \left\|\partial_{x_1}w\right\|_{BMO}
+\left\|\partial_{x_2}w\right\|_{BMO}\\
&\lesssim \left\|\frac{x_1}rw_r\right\|_{BMO}+
\left\|\frac{x_2}rw_r\right\|_{BMO}\\
&\lesssim \|w_r\|_{L^{\infty}}.
  \end{aligned}
\end{equation}
If we can find a such kind of slution, then possibly the modulation method (see the nice introduction in \cite{Lin20}) can be used to study the construction of the finite time blow up solutions to the Navier-Stokes equations.

\section{Acknowledgement}

The authors would like to thank Dr. Xiao Ren to point out that the solution of the Navier-Stokes equations corresponding to $V$ has no swirl.

The first author was partially supported by NSFC (12271487, 12171097), the second author was partially supported by NSFC (12171097), the Key Laboratory of Mathematics for Nonlinear Sciences (Fudan
University), Ministry of Education of China, Shanghai Key Laboratory for Contemporary Applied Mathematics, School of Mathematical Sciences, Fudan University, P.R. China,
and by Shanghai Science and Technology Program (21JC1400600).

\textbf{Data Availability Statement} Data sharing not applicable to this article as no
datasets were generated or analysed during the current study.

\textbf{Declarations}\\
\textbf{Conflict of interest} The authors have no relevant financial or non-financial interests
to disclose.

\bibliographystyle{plain}
\bibliography{reference}
	
\end{document}